\documentclass[amsppt,12pt]{amsart}
\usepackage{amsmath, amsthm, amssymb}
\usepackage[mathscr]{eucal}

\newcommand{\A}{\mathbb{A}}

\newcommand{\G}{\mathbb{G}}

\newcommand{\Q}{\mathbb{Q}}
\newcommand{\R}{\mathbb{R}}

\newcommand{\Z}{\mathbb{Z}}

\newcommand{\cF}{\mathcal{F}}

\newcommand{\scC}{\mathscr{C}}

\def\gg{\mathfrak{g}}

\newcommand{\gp}{\mathfrak{p}}

\newcommand{\gU}{\mathfrak{U}}
\newcommand{\gV}{\mathfrak{V}}
\newcommand{\gW}{\mathfrak{W}}

\newcommand{\BS}{\mathrm{B.S.}}

\newcommand{\cpct}{\mathrm{c}}

\newcommand{\Ext}{\mathrm{Ext}}

\newcommand{\la}{\mathrm{loc.an.}}

\newcommand{\limprojs}{\displaystyle \lim_{\stackrel{ \textstyle \leftarrow}{\scriptstyle s} }\ }
\newcommand{\limdr}{\displaystyle \lim_{\stackrel{ \textstyle \to}{\scriptstyle r} }\ }

\newcommand{\limpKp}{\displaystyle \lim_{\stackrel{\textstyle \leftarrow}{\scriptstyle K_{\gp} }\  }}
\newcommand{\limdKp}{\displaystyle \lim_{\stackrel{\textstyle \rightarrow}{\scriptstyle K_{\gp}}\ }}

\newcommand{\rank}{\mathrm{rank}}

\newcommand{\sing}{\mathrm{sing}}

\newcommand{\sq}{\hfill$\square$}

\newcommand{\Zps}{\Z/p^{s}}
\newcommand{\uZps}{\underline{\Zps}}

\newtheorem{theorem}{Theorem}

\newtheorem{corollary}{Corollary}

\theoremstyle{definition}

\theoremstyle{remark}

\begin{document}

\title{A reinterpretation of Emerton's $p$-adic Banach spaces}
\author{Richard Hill}
\date{14th July 2007}
\maketitle

\begin{abstract}
	It is shown that the $p$-adic Banach spaces introduced by
	Emerton are isomorphic to the cohomology groups of the sheaf of
	continuous $\Q_{p}$-valued functions on a certain space. Some
	applications of this result are discussed.
\end{abstract}

\section{Introduction}

Let $\G$ be a reductive group over a number field $k$.
We fix once and for all
 a maximal compact subgroup $K_{\infty}\subset \G(k\otimes\R)$,
 and we consider the ``Shimura manifolds'':
$$
	Y(K_{f})
	=
	\G(k) \backslash \G(\A) / K_{\infty}^{o}K_{f},
$$
where $K_{f}$ is a compact open subgroup of $\G(\A_{f})$,
and $K_{\infty}^{o}$ is the identity component in $K_{\infty}$.
Fix once and for all a finite prime $\gp$ of $k$, and let $p$ be the
rational prime below $\gp$.
In \cite{emerton} Emerton introduced the following spaces:
$$
	\tilde H^{\bullet}_{\ast}(K^{\gp},\Z_{p})
	=
	\limprojs
	\limdKp
	H^{\bullet}_{\ast}(Y(K_{\gp}K^{\gp}),\Zps),
$$
$$
	\tilde H^{\bullet}_{\ast}(K^{\gp},\Q_{p})
	=
	\tilde H^{\bullet}_{\ast}(K^{\gp},\Z_{p})\otimes_{\Z_{p}} \Q_{p},
$$
where $\ast$ is either the empty symbol, denoting usual cohomology,
or ``$\cpct$'', denoting compactly supported cohomology.
Here the $K_{\gp}$ ranges over the compact open subgroups
of $\G(k_{\gp})$, and $K^{\gp}$ is a fixed compact open subgroup of $\G(\A_{f}^{\gp})$.
The spaces $H^{\bullet}_{\ast}(K^{\gp},\Q_{p})$ are $p$-adic Banach spaces, and are
central to Emerton's construction of eigenvarieties in \cite{emerton}.

The aim of this paper is to give a more convenient interpretation
of these spaces.
To explain this interpretation, consider the topological space:
$$
	Y(K^{\gp})
	=
	\G(k)\backslash\G(\A)/K^{\gp}K_{\infty}^{o}
	=
	\limpKp
	Y(K^{\gp}K_{\gp}).
$$
Let $\scC_{\Z_{p}}$ (resp. $\scC_{\Q_{p}}$) be the sheaf of continuous
 $\Z_{p}$-valued (resp. $\Q_{p}$-valued) functions on $Y(K^{\gp})$.
The space $Y(K^{\gp})$ need not be compact, but it is homotopic to
a profinite simplicial complex, which we shall call $Y(K^{\gp})^{\BS}$.
We shall also use the notation
$\partial Y(K^{\gp})^{\BS}=Y(K^{\gp})^{\BS}\setminus Y(K^{\gp})$.

Our main result is the following.

\begin{theorem}
	There are canonical isomorphisms
	$$
		\tilde H^{\bullet}(K^{\gp},\Z_{p})
		=
		\check H^{\bullet}(Y(K^{\gp}),\scC_{\Z_{p}}),
		\quad
		\tilde H^{\bullet}_{\cpct}(K^{\gp},\Z_{p})
		=
		\check H^{\bullet}(Y(K^{\gp})^{\BS},\partial Y(K^{\gp}),\scC_{\Z_{p}}),
	$$
	and similarly for $\Q_{p}$.
	The right hand side of these equations is \v Cech cohomology, which in this case
	is equal to sheaf cohomology.
\end{theorem}

One of our aims in proving this result is to compare the spaces
 $\tilde H^{\bullet}$ and $\tilde H^{\bullet}_{\cpct}$.
An immediate consequence of our result is a long exact sequence
involving these two spaces
$$
	\to \tilde H^{n}_{\cpct}(K^{\gp},\Z_{p})
	\to \tilde H^{n}(K^{\gp},\Z_{p})
	\to \tilde H^{n}_{\partial} (K^{\gp},\Z_{p})
	\to \tilde H^{n+1}_{\cpct}(K^{\gp},\Z_{p})
	\to,
$$
where we are using the notation
$$
	\tilde H^{n}_{\partial} (K^{\gp},\Z_{p})
	=
	\check H^{n}(\partial Y(K^{\gp})^{\BS},\scC_{\Z_{p}}).
$$
This is significant, since one can show that
 $\tilde H^{n}_{\partial}$
vanishes unless $n$ is quite small.
For example, if $\G$ has real rank $1$
 and $\gp$ is the only prime of $k$ above $p$,
 then only $\tilde H^{0}_{\partial}$ is non-zero.
As a consequence, we know that for such groups,
 $\tilde H^{n}$ and $\tilde H^{n}_{\cpct}$ are equal for $n\ge 2$.
Generalizations of such results will be discussed in a forthcoming paper.

It is also envisioned that these results should give new insight into Eisenstein cohomology classes.
To see why this might be the case,
 we recall that Emerton proved a spectral sequence
$$
	\Ext^{p}_{\gg}(\check W, \tilde H^{q}_{\ast}(K^{\gp},k_{\gp})_{\la})
	\implies
	H^{p+q}_{\ast}(Y(K^{\gp}),W),
$$
where $W$ is a local system on $Y(K^{\gp})$ given by a finite dimensional representation
of $\G$ over $k_{\gp}$.
We remark that if $\G$ is semi-simple then
 $\tilde H^{n}_{\cpct}$ is zero for $n<\rank_{k}(\G)$;
 this follows easily by Poincar\'e duality.
On the other hand $\tilde H^{n}_{\partial}$ vanishes for all but small values of $n$.
Hence one might expect to be able to recover the boundary cohomology
quite low down in the filtration given by the spectral sequence.

\section{Some facts about \v Cech cohomology}

Let $X$ be a topological space and $\cF$ a presheaf on $X$.
For an open cover $\gU=\{U_{i}: i\in I\}$ of $X$,
 we define the \v Cech complex $\check C^{\bullet}(\gU,\cF)$ by
$$
	\check C^{n}(\gU,\cF)
	=
	\{(f_{i_{0},\ldots,i_{n}})_{i_{0},\ldots,i_{n}\in I^{n+1}} : 
	f_{i_{0},\ldots,i_{n}}\in \cF(U_{i_{0}}\cap \cdots \cap U_{i_{n}})\}.
$$
The cohomology groups of this complex are written $\check H^{\bullet}(\gU,\cF)$.
The \v Cech cohomology groups are defined to be
the direct limits of these cohomology groups:
$$
	\check H^{n}(X,\cF)
	=
	\lim_{\stackrel{\textstyle \to}{\scriptstyle \gU}} \check H^{n}(\gU,\cF).
$$
In fact $\check H^{n}(X,\cF)$ depends only on the sheafification of $\cF$.

\begin{theorem}[Leray's Theorem]
	Let $\cF$ be a sheaf on a topological space $X$
	and $\gU$ a countable open cover of $X$.
	If $\cF$ is acyclic on every finite intersection of
	elements of $\gU$, then
	$$
		\check H^{n}(X,\cF)
		=
		\check H^{n}(\gU,\cF).
	$$
\end{theorem}

\begin{theorem}[Thm. III.4.12 of \cite{bredon}]
	If $\cF$ is a sheaf on $X$ and $X$ is paracompact,
	 then the \v Cech cohomology groups of $\cF$ are equal to its sheaf
	cohomology groups, i.e. the derived functors of the global
	sections functor.
\end{theorem}

Given a presheaf $\cF$ on a topological space $Y$,
and a subspace $Z\subset Y$, we define presheafs $\cF_{Z}$ and $\cF^{Z}$
 on $X$ by
$$
	\cF_{Z}(U)
	=
	\begin{cases}
		\cF(U) & \hbox{if } U\cap Z\ne \emptyset\\
		0 & \hbox{otherwise.}
	\end{cases}
$$
$$
	\cF^{Z}(U)
	=
	\begin{cases}
		0 & \hbox{if } U\cap Z\ne \emptyset\\
		\cF(U) & \hbox{otherwise.}
	\end{cases}
$$
It turns out that $\check H^{\bullet}(Z,\cF)=\check H^{\bullet}(X,\cF_{Z})$,
and one defines
$$
	\check H^{\bullet}(Y,Z,\cF)
	=
	\check H^{\bullet}(Y,\cF^{Z}).
$$
There is a short exact sequence of presheafs:
$$
	0
	\to
	\cF^{A}
	\to
	\cF
	\to
	\cF_{A}
	\to
	0,
$$
This gives a long exact sequence:
$$
	\check H^{n}(\gU,\cF^{A}) \to
	\check H^{n}(\gU,\cF) \to
	\check H^{n}(\gU,\cF_{A}) \to
	\check H^{n+1}(\gU,\cF^{A}).
$$
Passing to the direct limit, we obtain
 the long exact sequence of \v Cech cohomology groups:
$$
	\check H^{n}(X,A,\cF) \to
	\check H^{n}(X,\cF) \to
	\check H^{n}(A,\cF) \to
	\check H^{n+1}(X,A,\cF).
$$
If $A$ is an abelian group, then we shall also write
$\underline A$ for the sheaf of locally constant $A$-valued functions.
Using Leray's theorem, one easily proves the following:

\begin{theorem}[Comparison Theorem]
	Let $Y$ be a finite simplicial complex,  and $Z\subset Y$
	a subcomplex.
	For any abelian group $A$,
	 we have
	$$
		\check H^{\bullet}(Y,Z, \underline A)
		=
		H^{\bullet}(Y,Z,A),
	$$
	where the right hand side is singular cohomology.
\end{theorem}

In fact the comparison theorem holds for much more general topological
spaces (see for example \cite{spanier}).

\begin{theorem}[Lem. 6.6.11, Cor 6.1.11 and Cor. 6.9.9 of \cite{spanier}]
	\label{spantheorem}
	Let $Y$ be a finite simplicial complex and
	$Z$ a subcomplex.
	For any abelian group $A$, we have
	$$
		H^{\bullet}_{\cpct}(Y\setminus Z,A)
		=
		H^{\bullet}(Y,Z,A).
	$$
\end{theorem}

\section{Proofs}

Emerton used the following formalism to introduce the
groups $\tilde H^{n}_{\ast}$.
Let $G$ be a compact, $\Q_{p}$-analytic group,
and fix a basis of open, normal subgroups:
$$
	G
	=
	G_{0}
	\supset
	G_{1}
	\supset
	\ldots.
$$
Suppose we have a sequence of simplicial maps between
finite simplicial complexes
$$
	\cdots \to Y_{2}
	\to Y_{1}
	\to Y_{0},
$$
and subcomplexes:
$$
	\cdots \to Z_{2} \to Z_{1} \to Z_{0},
$$
each equipped with a right action of $G$, and satisfying the following
conditions:
\begin{enumerate}
	\item
	the maps in the sequence are $G$-equivariant;
	\item
	$G_{r}$ acts trivially on $Y_{r}$.
	\item
	if $0\le r' \le r$ then the maps $Y_{r}\to Y_{r'}$ and $Z_{r}\to Z_{r'}$
	are Galois covering maps with deck transformations
	provided by the natural action of $G_{r'}/G_{r}$ on $Y_{r}$.
\end{enumerate}
Given this data, we let $Y$ be the projective limit of the spaces $Y_{r}$,
 and $Z$ be the projective limit of the spaces $Z_{r}$.
We shall use the notation $Y^{0}=Y\setminus Z$,
 $Y_{i}^{0}=Y_{i}\setminus Z_{i}$.
Emerton defined the following spaces:
$$
	\tilde H^{n}_{\cpct}(Y^{0},\Z_{p})
	=
	\limprojs
	\limdr
	H^{n}_{\cpct}(Y_{r}^{0},\Zps),
	\quad
	\tilde H^{n}_{\cpct}(Y^{0},\Q_{p})
	=
	\tilde H^{n}_{\cpct}(Y^{0},\Z_{p})
	\otimes_{\Z_{p}}
	\Q_{p}.
$$
In applications, $G$ will be a compact open subgroup of $\G(k_{\gp})$,
and the space $Y$ will be either $Y(K^{\gp})^{\BS}$
or $\partial Y^{\BS}$. If $Y=Y(K^{\gp})^{\BS}$ then we may use the subspace
$Z=\partial Y(K^{\gp})^{\BS}$.

\begin{theorem}
	With the above notation,
	$$
		\tilde H^{n}_{\cpct}(Y^{0},\Z_{p})
		=
		\check H^{n}(Y, Z, \scC_{\Z_{p}}),
		\quad
		\tilde H^{n}_{\cpct}(Y^{0},\Q_{p})
		=
		\check H^{n}(Y,Z, \scC_{\Q_{p}}).
	$$
\end{theorem}

\proof
We shall prove the case with coefficients in $\Z_{p}$.
The $\Q_{p}$ case is a consequence.
We shall write $\scC$ instead of $\scC_{\Z_{p}}$.
To prove the theorem,
 we construct an acyclic cover of $Y$ and apply Leray's Theorem.

\subsection{A cover}

We first choose a finite open cover $\gU$
 of $Y_{0}$ with the following properties:
\begin{enumerate}
	\item
	If $U$ is an intersection of finitely many sets in $\gU$
	then either $U$ is empty or $U$ is contractible.
	\item
	If $U$ is an intersection of finitely many sets in $\gU$
	and $U\cap Z_{0}$ is non-empty,
	then $U\cap Z_{0}$ is a deformation retract of $U$.
\end{enumerate}
For each $U\in \gU$, we let $U^{(r)}$ be the preimage
of $U$ in $Y_{r}$.
The sets $U^{(r)}$ form an open cover $\gU^{(r)}$ of $Y_{r}$,
and have the following properties:
\begin{enumerate}
	\item
	For every $U_{1}^{(r)},\ldots,U_{s}^{(r)}\in \gU^{(r)}$
	with non-empty intersection, the intersection
	$U_{1}^{(r)}\cap \cdots \cap U_{s}^{(r)}$
	is isomorphic as a topological $G$-set to
	$(U_{1}\cap \cdots \cap U_{r})\times (G/G_{r})$.
	In particular, the intersection is homotopic to a
	finite set.
	\item
	If $U_{1}^{(r)},\ldots,U_{s}^{(r)}\in \gU^{(r)}$
	and $U_{1}^{(r)}\cap \cdots \cap U_{s}^{(r)}\cap Z_{r}$ is
	non-empty,
	then	$U_{1}^{(r)}\cap \cdots \cap U_{s}^{(r)}\cap Z_{r}$
	is a deformation retract of $U_{1}^{(r)}\cap \cdots \cap U_{s}^{(r)}$.
\end{enumerate}
Furthermore, for each set $U\in \gU$, we define $\tilde U$
to be the preimage of $U$ in $Y$.
The sets $\tilde U$ form an open cover $\tilde\gU$ of $Y$.
We immediately verify the following:
\begin{enumerate}
	\item
	if $\tilde U_{1},\ldots,\tilde U_{s}\in\tilde \gU$ have non-empty intersection,
	then their intersection is equivalent as a topological $G$-set
	to $(U_{1}\cap \cdots\cap U_{s})\times G$.
	\item
	if $\tilde U_{1},\ldots,\tilde U_{s}\in\tilde \gU$ and
	$\tilde U_{1} \cap\cdots \cap\tilde U_{s}\cap Z$ is non-empty,
	then $\tilde U_{1} \cap\cdots \cap\tilde U_{s}\cap Z$ is a deformation retract
	of $\tilde U_{1} \cap\cdots \cap\tilde U_{s}$.
\end{enumerate}

\subsection{$\gU^{(r)}$ is $(\uZps)^{Z_{r}}$-acyclic}

Let $U$ be an intersection of finitely many sets in $\gU$,
and let $U^{(r)}$ be the preimage of $U$ in $Y_{r}$.
We know that $U$ is contractible, and $U^{(r)}=U\times (G/G_{r})$.
The sheaf $(\uZps)^{Z_{r}}$ on $Y_{r}$ consists of locally constant
$\Zps$-valued functions, which vanish on $Z_{r}$.
It follows that $\check H^{\bullet}(U^{(r)},(\uZps)^{Z_{r}})$ is a direct sum of finitely
many copies of $\check H^{\bullet}(U,(\uZps)^{Z_{0}})$.
We must therefore show that $\check H^{n}(U,(\uZps)^{Z_{0}})=0$ for all $n>0$.

If $U\cap Z_{0}$ is empty, then we have
$\check H^{n}(U,(\uZps)^{Z_{0}})=\check H^{n}(U,\uZps)$.
By the comparison theorem, this is the same as singular cohomology, and therefore
only depends on $U$ up to homotopy.
Since $U$ is contractible, it follows that $\check H^{n}(U,\uZps)=0$ for $n>0$.

Suppose instead that $U\cap Z_{0}$ is non-empty.
In this case, we know that $U\cap Z_{0}$ is a deformation retract of $U$.
It follows that the restriction map
$H^{\bullet}_{\sing}(U,\Zps) \to H^{\bullet}_{\sing}(U\cap Z_{0},\Zps)$ is an isomorphism.
By the comparison theorem, it follows that the map
$\check H^{\bullet}(U,\uZps)
 \to \check H^{\bullet}(U\cap Z_{0},\uZps)$ is an isomorphism.
The long exact sequence shows that $\check H^{\bullet}(U,U\cap Z_{0},\uZps)=0$.

In particular, using Leray's Theorem, we have
$$
	\check H^{\bullet}(Y_{r},Z_{r},\uZps)
	=
	\check H^{\bullet}(\gU^{(r)},(\uZps)^{Z_{r}}).
$$

\subsection{$\tilde \gU$ is $\scC$-acyclic}

Let $U$ be an intersection of finitely many sets in $\gU$,
and let $\tilde U$ be the preimage of $U$ in $Y$.
We know that $U$ is contractible, and $\tilde U=U\times G$.
We must show that $\check H^{n}(\tilde U,\scC)=0$ for $n>0$.

Let $\tilde \gV$ be an open cover of $\tilde U$, and choose an element
 $\sigma\in \check H^{n}(\tilde \gV,\scC)$ with $n>0$.
We shall find a refinement $\tilde \gW$ of $\tilde \gV$, such that
 the image of $\sigma$ in $\check H^{\bullet}(\tilde \gW,\scC)$ is zero.
By passing to a refinement of $\tilde \gV$ if necessary, we may assume that
$\tilde\gV$ is finite, and that each element of $\tilde \gV$
 is of the form $V_{i}\times H_{i}$
 for some open subset $V_{i}\subset U$ and some open coset $H_{i}\subset G$.
By refining still further, we may assume that the cosets $H_{i}$ are all cosets of
the same open subgroup $G_{r}\subset G$.
This means that $\tilde \gV$ is the pullback
 of an open cover $\gV^{(r)}$ of $U^{(r)}$.
For an open subset $V^{(r)}\subset U^{(r)}$, we have
$$
	\scC(\tilde V^{(r)})
	=
	S(V^{(r)}),
$$
where $S$ is the locally constant sheaf on $U^{(r)}$
with values in $\scC(G_{r},\Z_{p})$
and $\tilde V^{(r)}$ is the preimage of $V^{(r)}$ in $\tilde U$.
It follows that
$$
	\check H^{\bullet}(\tilde\gV, \scC_{\Z_{p}})
	=
	\check H^{\bullet}(\gV^{(r)}, S).
$$
Since $U^{(r)}$ is homotopic to a finite set and $S$ is a constant sheaf,
 it follows that $\check H^{>0}(U^{(r)}, S)$ is zero.
This implies there is a refinement $\gW^{(r)}$ of $\gV^{(r)}$, such that
 the image of $\sigma$ in $\check H^{\bullet}(\gW^{(r)}, S)$ is zero.
Pulling $\gW^{(r)}$ back to $\tilde U$,
 we have a refinement $\tilde \gW$ of $\tilde\gV$,
 such that the image of $\sigma$ in $\check H^{\bullet}(\tilde\gW,\scC_{\Z_{p}})$
 is zero.

\subsection{$\tilde \gU$ is $\scC^{Z}$-acyclic}

Let $U$ be an intersection of finitely many sets in $\gU$,
and let $\tilde U$ be the preimage of $U$ in $Y$.
We know that $U$ is contractible, and $\tilde U=U\times G$.
We must show that $\check H^{n}(\tilde U,\tilde U\cap Z,\scC)=0$ for $n>0$.
If $U$ does not intersect $Z_{0}$,
 then this follows from the previous part of the proof.
We therefore assume that $U$ intersects $Z_{0}$.
In this case, we know that $U\cap Z_{0}$ is a deformation retract of $U$.
In particular, $U\cap Z_{0}$ is contractible, and
$\tilde U\cap Z= (U\cap Z_{0})\times G$.
The previous part of the proof shows that $\check H^{>0}(\tilde U, \scC)=0$
and $\check H^{>0}(\tilde U \cap Z,\scC)=0$.
Furthermore, one sees immediately that the restriction map
 $\check H^{0}(\tilde U,\scC)\to \check H^{0}(\tilde U \cap Z,\scC)$
 is an isomorphism.
Hence by the long exact sequence,
 we have $H^{\bullet}(\tilde U,\tilde U \cap Z,\scC)=0$.

Thus by Leray's Theorem, we have:
$$
	\check H^{\bullet}(Y,Z,\scC)
	=
	\check H^{\bullet}(\tilde \gU,\scC^{Z}).
$$

\subsection{}
Fix for a moment a cohomological degree $n$,
and let $U_{1},\ldots,U_{N}$ be the non-empty intersections
of $n+1$-tuples of sets in $\gU$, for which $U_{i}\cap Z_{0}=\emptyset$.
For each $U_{i}$, we let $U_{i}^{(r)}$ be the preimage of $U_{i}$ in $Y_{r}$
and $\tilde U_{i}$ be the preimage of $U_{i}$ in $Y$.

Recall that $\check H^{\bullet}(\gU^{(r)},(\Z/p^{r})^{Z})$
 is the cohomology of the chain complex
$$
	\check C^{n}(\gU^{(r)},(\uZps)^{Z_{r}})
	=
	\prod_{i=1}^{N} (\uZps)^{Z_{r}}(U_{i}^{(r)}),
$$
Each $U_{i}$ is contractible and disjoint from $Z_{0}$.
Furthermore $U_{i}^{(r)}=U_{i}\times (G/G_{r})$,
 so we have an isomorphism of $G$-modules:
$(\uZps)^{Z_{r}}(U_{i}^{(r)})=(\uZps)(G/G_{r})$.
This gives
$$
	\check C^{n}(\gU^{(r)},(\uZps)^{Z})
	=
	(\uZps)(G/G_{r})^{N},
$$
Similarly, we have 
$$
	\check C^{n}(\tilde \gU,(\uZps)^{Z})
	=
	\left((\uZps)(G)\right)^{N}.
$$
Comparing the two formulae, it is clear that
$$
	\check C^{\bullet}(\tilde \gU,(\uZps)^{Z})
	=
	\lim_{\to\atop r} \check C^{\bullet}(\gU^{(r)},(\uZps)^{Z_{r}}).
$$
Since the functor $\limdr$ is exact, we have
$$
	\check H^{\bullet}(\tilde \gU,(\uZps)^{Z})
	=
	\lim_{\to\atop r} \check H^{\bullet}(\gU^{(r)},(\uZps)^{Z_{r}}).
$$

\subsection{}
Note also that $\scC^{Z}(\tilde U_{i})=\scC(G)$,
and so we have
$$
	\check C^{n}(\tilde \gU,\scC^{Z})
	=
	\scC(G)^{N}.
$$
It follows that $\check C^{n}(\tilde \gU,\scC^{Z})$
 is an admissible $\Z_{p}[G]$-module in the sense of \cite{emerton}.
Furthermore we have:
$$
	\check C^{\bullet}(\tilde \gU,\scC^{Z})
	=
	\limprojs \check C^{\bullet}(\tilde \gU,(\uZps)^{Z}),
	\quad
	\check C^{\bullet}(\tilde \gU,(\uZps)^{Z})
	=
	\check C^{\bullet}(\tilde \gU,\scC^{Z})/p^{s}.
$$
Hence by Proposition 1.2.12 of \cite{emerton}, we have:
$$
	\check H^{\bullet}(\tilde \gU,\scC^{Z})
	=
	\limprojs \check H^{\bullet}( \tilde \gU,(\uZps)^{Z}).
$$
By the previous part of the proof, we have:
$$
	\check H^{\bullet}(\tilde \gU,\scC^{Z})
	=
	\limprojs \limdr
	\check H^{\bullet}( \gU^{(r)},(\uZps)^{Z_{r}}).
$$
Since our covers are acyclic, this translates to
$$
	\check H^{\bullet}(Y,Z,\scC)
	=
	\limprojs \limdr
	\check H^{\bullet}( Y_{r},Z_{r},\uZps).
$$
On the other hand, by Theorem \ref{spantheorem},
we have
$$
	\check H^{\bullet}( Y_{r},Z_{r},\uZps)
	=
	H^{\bullet}_{\cpct}( Y_{r}\setminus Z_{r},\Zps).
$$
The result follows.
\sq

\begin{corollary}
	With the above notation,
	$$
		\tilde H^{n}(Y,\Z_{p})
		=
		\check H^{n}(Y, \scC_{\Z_{p}}),
		\quad
		\tilde H^{n}(Y,\Q_{p})
		=
		\check H^{n}(Y, \scC_{\Q_{p}}).
	$$
\end{corollary}

\proof
We apply the theorem in the case that $Z$ is empty.
Since $Y^{0}=Y$, which is compact,
 it follows that usual cohomology is the same as compactly supported cohomology on
 each $Y_{r}$.
\sq

\begin{corollary}
	In the notation of the introduction, there are long exact sequences:
	$$
		\tilde H^{n}_{\cpct}(K^{\gp},\Z_{p})
		\to
		\tilde H^{n}(K^{\gp},\Z_{p})
		\to
		\tilde H^{n}_{\partial}(K^{\gp},\Z_{p})
		\to
		\tilde H^{n+1}_{\cpct}(K^{\gp},\Z_{p}),
	$$
	$$
		\tilde H^{n}_{\cpct}(K^{\gp},\Q_{p})
		\to
		\tilde H^{n}(K^{\gp},\Q_{p})
		\to
		\tilde H^{n}_{\partial}(K^{\gp},\Q_{p})
		\to
		\tilde H^{n+1}_{\cpct}(K^{\gp},\Q_{p}).
	$$
\end{corollary}

\proof
For convenience, we shall write $Y$ instead of $Y(K^{\gp})$.
We have shown above that
\begin{eqnarray*}
	\tilde H^{\bullet}(K^{\gp},\Q_{p})
	&=&
	\check H^{\bullet}(Y,\scC_{\Q_{p}}),\\
	\tilde H^{\bullet}_{\partial}(K^{\gp},\Q_{p})
	&=&
	\check H^{\bullet}(\partial Y^{\BS},\scC_{\Q_{p}}),\\
	\tilde H^{\bullet}_{\cpct}(K^{\gp},\Q_{p})
	&=&
	\check H^{\bullet}(Y^{\BS},\partial Y^{\BS},\scC_{\Q_{p}}).
\end{eqnarray*}
There is a long exact sequence in \v Cech cohomology:
$$
	\check H^{n}(Y^{\BS},\scC_{\Q_{p}})
	\to
	\check H^{n}(\partial Y,\scC_{\Q_{p}})
	\to
	\check H^{n+1}(Y^{\BS},\partial Y^{\BS},\scC_{\Q_{p}})
	\to
	\check H^{n+1}(Y^{\BS},\scC_{\Q_{p}}).
$$
The same holds for $\Z_{p}$.
\sq


\begin{thebibliography}{10}

\bibitem{bredon}
G.~E.~Bredon.
\newblock {\em Sheaf Theory, second edition}.
\newblock Graduate Texts in Mathematics vol. 170.
\newblock Springer-Verlag 1997.

\bibitem{emerton}
M.~Emerton.
\newblock On the interpolation of systems of eigenvalues attached to
  automorphic {H}ecke eigenforms.
\newblock {\em Invent. Math.}, 164(1):1--84, 2006.



\bibitem{spanier}
E.~H. Spanier.
\newblock {\em Algebraic Topology}.
\newblock Springer-Verlag New York, inc., 1966.


\end{thebibliography}
\end{document}